\documentclass{amsart}
\allowdisplaybreaks[4]

\newcommand{\R}{\mathbf{R}}

\newcommand{\sig}{\sigma}

\newcommand{{\ba}}{\bf a}
\newcommand{\ve}{\varepsilon}

\newcommand{\ga}{\gamma}
\newcommand{\pa}{\partial}
\newcommand{\ra}{\rightarrow}
\newcommand{\Om}{\Omega}
\newcommand{\del}{\delta}

\newcommand{\De}{\Delta}
\newcommand{\al}{\alpha}

\newcommand{\be}{\begin{equation}}
\newcommand{\ee}{\end{equation}}

\newcommand{\nn}{\nonumber}

\newtheorem{lem}{Lemma}{\bf}{\it}
{\it}{\rm}

\newtheorem{theorem}{Theorem}
\newtheorem{proposition}{Proposition}

\numberwithin{theorem}{section}
\numberwithin{lem}{section}
\numberwithin{equation}{section}
\numberwithin{proposition}{section}
\numberwithin{corollary}{section}
\title[Efficient Operation of a Bank]
{On a Model for the Efficient Operation of a Bank or Insurance Company}

\author{Joseph G. Conlon and Hyekyung Min}

\address{University of Michigan\\ Department of Mathematics\\ Ann Arbor,
  MI 48109-1109}
\email{conlon@umich.edu, hyekyung@umich.edu}


\begin{document}

\maketitle

\begin{abstract}
In this paper the authors study a model for the optimal operation of a bank or insurance company which was recently introduced by Peura and Keppo. The model generalizes a previous one of Milne and Robertson by allowing the bank  to raise capital as well as to pay out dividends. Optimal operation of the bank is determined by solving an optimal control problem. In this paper it is shown that the solution of the optimal control problem proposed by Peura and Keppo exists for all values of the parameters and is unique.

\end{abstract}

\section{Introduction}

\ \ \ \ \ In this paper we study a model for the optimal operation of a bank or insurance company which was introduced by Peura and Keppo \cite{pk}.  In this model capital is invested in a risky asset whose evolution is described by Brownian motion with drift.  Thus if $X(t)$ is the bank's capital at time $t$, then
\be \label{A1}
dX(t) = \mu dt + \sig dW(t),
\ee
where $W(t)$ is Brownian motion, $\mu > 0$ is the drift and $\sig > 0$ the volatility.  For an insurance company model $\mu$ represents the expected premium collection rate minus the expected claims payment rate.

In addition to investing in the asset described by (\ref{A1}) the bank also pays dividends to its owners and raises capital from them.  Dividend payments can be implemented instantaneously, but capital issuance is associated with a delay of length $\De$ and a fixed cost $K$.  If capital is ordered at time $t$ it is actually received at time $t+\De$.  Since there is a fixed cost $K$ associated with ordering the capital, if an amount $s$ of capital is ordered  then $s-K$ of that goes to increasing the total capital of the bank.  Furthermore, the model allows the owners to decide on the amount of capital based on all information up to time $t+\De$. The fixed cost $K$ is also paid at time $t+\De$.

The payment of dividends and the raising of capital is controlled by a policy $\pi$.  For $t > 0$ let $L^\pi(t)$ be defined by
\be \label{B1}
L^\pi(t) = \ {\rm dividends \; paid \; out \ - \ capital \; raised \; up \; to\; time} \ t .
\ee
It is assumed that capital is raised at a set of discrete times $t^\pi_1 < t^\pi_2 < t^\pi_3 < \dots$, where the number $N^\pi(t)$ of times capital is raised up to time $t$ is given by the formula,
\be \label{C1}
N^\pi(t) = \sup \{ i : t^\pi_i \le t \}.
\ee
It is also assumed that the times between raising capital is at least $\De$ and that no dividends are paid out during the period $\De$ before a capital raising time.  Thus, assuming $L^\pi(t)$ is a right continuous process, $L^\pi(t)$ is constant in the interval $\big( t^\pi_i - \De, \; t^\pi_i \big)$. At the time $t^\pi_i$ the bank decides the amount of capital it needs to raise. If its portfolio has performed particularly well in the previous time period of length $\De$  then it may not raise capital but actually pay out dividends at time  $t^\pi_i$. In all cases it has to pay the capital raising cost $K$.  Since dividend payments are non-negative the process $L^\pi(t)$ is increasing for $t \not\in \{ t^\pi_i : i = 1,2,...\}$.  

Let $X^\pi(t)$ be the amount of capital the bank has at time $t$.  Then from (\ref{A1}), (\ref{B1}), (\ref{C1}) the evolution of $X^\pi(t)$ is governed by the equation,
\be \label{E1}
d X^\pi(t) = \mu dt + \sig dW(t) - dL^\pi(t) - Kd N^\pi(t).
\ee
As is usual in control theory it is assumed that the policy $\pi$ depends only on information up to the present time.  Hence one assumes that the process $L^\pi(t)$ is right continuous and measurable with respect to the $\sig$ field generated by $X^\pi(s), \; s < t$.

The main concern of \cite{pk} is to determine a policy $\pi$ which maximizes the expected payments to the owners of the bank over its life time.  They therefore define a value function $V^\pi(x)$, $x\ge 0$, by
\be \label{F1}
V^\pi(x) = E\left[ \int^{\tau_\pi}_0 e^{-\rho t} dL^\pi(t)  \Big | X^\pi(0) = x \right],
\ee
where $\tau_\pi = \sup \left\{ t > 0 : X^\pi(t) > 0 \right\}$ is the life time of the bank.  The optimal value function is then given by
\be \label{G1}
V(x) = \sup_\pi \; V^\pi(x),
\ee
where the supremum norm is taken over all allowable strategies.  In \cite{pk} an expression for the function $V$ of (\ref{G1}) is obtained and a corresponding strategy to realize it.  The function $V$ is characterized by 2 parameters $u_{1,K}(\De), \; u_{2,K} (\De)$ satisfying $0 \le u_{1,K}(\De) < u_{2,K}(\De) < \infty$.  If the capital the bank holds is less than $u_{1,K}(\De)$ then a capital raising event is initiated.  If the capital exceeds $u_{2,K}(\De)$ then dividends are paid out.  If the capital the bank holds lies between $u_{1,K}(\De)$ and $u_{2,K} (\De)$ then it is fully invested in the risky asset described by (\ref{A1}).  A notable feature of the model is that $\tau_\pi < \infty$ with probability 1 for an optimal strategy.  Thus maximally efficient operation of the bank gives rise to zero capital in finite time.

The model of Peura and Keppo generalizes an earlier model of Milne and Robertson \cite{mr} which allows dividend payments but not raising of capital (see also \cite{ht,ht1}).    The Milne-Robertson model can be recaptured from that of Peura and Keppo by simply taking the cost $K$ of raising capital to be sufficient large.  In that case $u_{1,K}(\De) = 0$ and $u_{2,K} (\De) = u_0$, where $u_0$ is the Milne-Robertson threshold for the payment of dividends.  The value function $V$ of (\ref{G1}) now satisfies the equation
\be \label{H1}
(A-\rho) V(x)=0, \; 0 < x < u_0, \; V(0) = 0, \; V'(u_0) = 1,
\ee
where $A$ is the infinitesimal generator for the process (\ref{A1}),
\be \label{I1}
A = \frac 1 2 \sig^2 \frac{d^2}{dx^2} + \mu \ \frac d{dx},
\ee
and $u_0 > 0$ is chosen so that the function $V$, when continued by a linear function for $x > u_0$, is $C^2$.  The complete value function is then this linearly extended function which satisfies (\ref{H1}) for $0 \le x \le u_0$.

When $K$ and $\De$ are small enough then $u_{1,K}(\De) >0$ and it is determined from the solution of a free boundary problem which is a zero latent heat limit for the Stefan problem (see \cite{fs} for a Stefan problem occuring in finance).  The approach in \cite{pk} to obtaining the solution to (\ref{G1}) is to study the Bellman inequalities \cite{fr, fls} corresponding to the optimal control problem.  These are given as follows:
\begin{itemize} 
 \item[(a)]  \    $V(0) = 0, $
\item[(b)]  \  $V(x) \ge MV(x), $
\item[(c)]   \    $(A - \rho) V(x) \le 0, $ \be \label{J1} \ee
\item[(d)]  \ $V'(x) \ge 1, $
\item[(e)]   \  $[V(x) - MV(x)][(A-\rho) V(x)][V'(x) - 1]=0,$
\end{itemize}
where the operator $M$ in (b) is defined by
\be \label{K1}
Mf(x) = E\left[ e^{-\rho \De} \sup_s \Big[f\Big(X(\De) + s \Big) - s - K \Big] I_{\tau > \De} | X(0) = x \right], \ x > 0.
 \ee
Here $X(t)$ is the diffusion (\ref{A1}) and $\tau$ is the first hitting time at 0.  Our main theorem is that the Peura-Keppo solution to (\ref{G1}) is the unique solution to the system of inequalities  (\ref{J1}).
\begin{theorem}
 For all $K,\De > 0$ there is a unique $C^1$ solution $V(x)$ to the system of inequalities (\ref{J1}).  Further, for $\ve > 0$ there is a control policy $\pi = \pi_\ve$ such that if $V_\ve$ is defined by (\ref{F1}) then $\lim_{\ve\ra 0} V_\ve (x) = V(x), \ x \ge 0$.
\end{theorem}

In $\S 2$ we construct the solution to (\ref{J1}).  The main mathematical fact needed is that a solution $u(x,t), \; t>0$, of the diffusion equation such that the initial data $u(x,0)$ has just one change of sign, has at most one sign change for all $t> 0$.  This property of the diffusion equation has been used previously \cite{b, mck}, and a proof based on probability given.  In the appendix we give a proof using the maximum principle \cite{pw}.  In $\S 3$ we give the proof of uniqueness and construct the policies $\pi_\ve$.  Here we need also to use the fact that solutions to $u(x,t)=0$ are non degenerate.

The final section of the paper is concerned with studying the thresholds $u_{1,K}(\De)$ and $u_{2,K}(\De)$, in particular their asymptotic behavior as $\De \ra 0$.  We also show that if $u_{1,K}(\De) > 0$ the optimal function $V(x)$, which is $C^\infty$ for $x \not= u_{1,K}(\De)$, $u_{2,K}(\De)$ and $C^2$ at $x = u_{2,K}(\De)$, is not $C^2$ at $u_{1,K}(\De)$.

\section{Solution to System of Bellman Equations}
\ \ In this section we construct a solution to the system of inequalities (\ref{J1}).  To do this we first consider solutions $V(x)$ to the equation
\be \label{A2}
(A-\rho) V(x) = 0, \ \ x \in \R.
\ee
For any $u_0 > 0$, there is a unique solution $V_0$ to (\ref{A2}) with the initial conditions,
\be \label{B2}
V_0(u_0) = \mu / \rho, \ \ \ V'_0(u_0) = 1.
\ee
Evidently (\ref{A2}), (\ref{B2}) imply that $V_0''(u_0) = 0$.  Since solutions to (\ref{A2}) have just one point of inflection it follows that the function $V_0$ is concave for $x <u_0$ and convex for $x > u_0$.  Evidently by translation there is a unique $u_0$ such that $V_0(0) = 0$.  This unique $u_0$ is given in \cite{pk} by the formula,
\be \label{C2}
u_0 = \frac 1{r_1 + r_2} \ \ell n \left[ \frac{\rho + \mu r_2}{\rho - \mu r_1}\right],
\ee
where $r_1, r_2$ are the characteristic roots for equation (\ref{A2}),
\[	r_1 = \frac{-\mu + \sqrt{\mu^2 + 2\sig^2 \rho}}{\sig^2}, \ \ 
        r_2 = \frac{\mu + \sqrt{\mu^2 + 2\sig^2 \rho}}{\sig^2} .
\]
Note that since $V_0$ is concave for $x < u_0$ and $V_0'(u_0) = 1$ it follows that $u_0$ of (\ref{C2}) satisfies the inequality,
\be \label{D2} 
0 < u_0 < \mu/\rho.
\ee
Next for the diffusion process with generator $A$, started at $x > 0$, let $\tau_x$ be the first hitting time at 0.  We define for $x > 0, \ t > 0$, the function $p(x,t)$ by
\be \label{E2}
p(x,t) = P(\tau_x < t).
\ee
Then $p(x,t)$ is a solution to the equation,
\be \label{F2}
\frac{\pa p}{\pa t}=A p,\quad x>0, \ t>0,
\ee
with boundary and initial conditions,
\be \label{G2}
p(0,t)=1, \ t>0; \quad p(x,0)=0, \ x>0.
\ee
It is also evident from the representation (\ref{E2}) that $p(x,t)$ satisfies the inequalities,
\be \label{H2}
\frac {\pa p}{\pa t} \ge 0, \ \ \frac {\pa p}{\pa x} \le 0, \ \ x > 0, \ t > 0.
\ee
It follows from (\ref{F2}), (\ref{H2}) that for any fixed $t > 0$ the function $p(x,t)$ is a convex function of $x, x>0$.

Let $h(x,t),  \  x > 0, t > 0$, be a solution of the equation,
\[	\frac{\pa h}{\pa t} = Ah, \ \ \ x >0, \ t >0,   \]
with boundary and initial conditions
\[	h(0,t) = 0, \ t > 0; \ \ h(x,0) = x, \ x > 0.	\]
It is easy to see that $h$ and $p$ are related by the formula,
\be \label{K2}
h(x,t) = x + \mu t - \mu \int^t_0 \; p(x,s)ds.
\ee
\begin{lem}
For fixed $t>0$ the function $h(x,t)$ is a concave function of $x, \ x > 0$.  It also satisfies the inequalities,
\be \label{I2}
\mu/\rho - u_0 + h(x,t) \le e^{\rho t} [x - u_0 + \mu/\rho], \ x \ge u_0, 
\ee
\be \label{J2}
h(x,t) \le e^{\rho t} V_0(x), \ x \ge 0.
\ee
\end{lem}
\begin{proof}   The concavity of $h$ follows from (\ref{K2}) and the convexity of $p$.  The inequality (\ref{I2}) follows from (\ref{K2}).  The inequality  (\ref{J2}) follows from the maximum principle for the diffusion equation since $V_0(x) \ge x, \ x>0$. 
\end{proof}
\begin{proposition}
Suppose $u_0$ is given by (\ref{C2}) and the cost $K$ of capital issuance satisfies $K \ge \mu/\rho - u_0$.  Define $V(x)$ by
\be \label{L2}
V(x) = V_0(x), \ 0 \le x \le u_0; \ V(x) = x-u_0 +\mu /\rho, \; x > u_0.
\ee
Then $V(x)$ is a $C^2$ function and satisfies the system of inequalities (\ref{J1}).
\end{proposition}
\begin{proof}  It is easy to see that the function $V(x)$ of (\ref{L2}) satisfies (\ref{J1}) (a), (c), (d), (e).  To prove (b) we note that the function $MV(x)$ of (\ref {K1}) is given by the expression,
\be \label{P2}
MV(x) = e^{-\rho \De} \left\{ \left[ \frac \mu \rho - u_0 - K\right] \{1-p(x,\De)\} + h(x,\De)\right\}.
\ee
Now (b) follows from Lemma 2.1 on noting (\ref{D2}). 
\end{proof}
From (\ref{J2}) we have that
\[	\frac{\pa h}{\pa x}(0,t) \le e^{\rho t} V'_0(0), \ t > 0.	\]
Observe also that by the Hopf maximum principle \cite{pw} one has $\pa p/\pa x(0,t) < 0$, $t >0$.  Hence one may define a smooth function $\beta(t), \; t>0$, by the formula,
\be \label{M2}
\beta(t) = \left[ \frac{\pa h}{\pa x}(0,t) - e^{\rho t} V'_0(0)\right] \Big/ \frac{\pa p}{\pa x}(0,t) .
\ee
\begin{lem}
The function $\beta(t)$ of (\ref{M2}) is strictly monotonic increasing and satisfies
\[	\lim_{t \ra 0} \beta(t) = 0, \ \ \ \lim_{t \ra \infty} \beta(t) = \infty. \]
\end{lem}
\begin{proof}  Observe by the Hopf maximum principle that $\beta(t) > 0,  \ t > 0$.  To prove monotonicity let $T>0$ and consider the function $u(x,t)$ defined by
\[	u(x,t) = \beta(T) \{1-p(x,t)\} + h(x,t) - e^{\rho t} V_0(x).	\]
Then $u(x,t)$ is a solution of the diffusion equation (\ref{F2}) with initial and boundary conditions,
\[	u(x,0)= \beta(T) + x - V_0(x), \ \ \ u(0,t) = 0.	\]
Evidently $u(x,0)$ is a monotonically decreasing function satisfying,
\[	\lim_{x\ra 0} u(x,0) =  \beta(T), \ \ \ \lim_{x \ra \infty} u(x,0) = -\infty. \]
In particular $u(x,0)$ has exactly one sign change.  It follows therefore by Theorem A1 of the appendix that $u(x,t)$ has at most one sign change for any fixed  $t > 0$.

Consider now the function $u(x,T)$.  Suppose first that $u(x,T) \le 0, \ x \ge 0$.  Then by the Hopf maximum principle  $u(x,t) < 0, \; x > 0$, for any $t > T$ and $\pa u / \pa x(0,t) < 0$.   Thus $\beta(t) > \beta(T)$.  Alternatively there is an interval $(0,\al)$ for which $u(x,T) > 0, \; x \in (0,\al)$.  By the Hopf principle one must have then $\pa u/ \pa x (0,T) > 0$ which contradicts the definition of $\beta(T)$.  We have show that $\beta(t)$ is strictly monotonic increasing.

To find the limit of $\beta(t)$ as $t \ra 0$ we compute the limits of the numerator and denominator of (\ref{M2}).  For the numerator we clearly have that
\be \label{N2}
\lim_{t \ra 0} \left[ \frac{\pa h}{\pa x}(0,t) - e^{\rho t} V'_0(0)\right] = 1 - V'_0(0) < 0.	
\ee
To find the limit of the denominator we use the Green's functions $G(x,y,t)$ for the equation (\ref{F2}).  Thus by the reflection principle we have that
\begin{eqnarray*}
G(x,y,t) &=& \frac 1{\sqrt{2\pi \sig^2\;t}} \Bigg\{ \exp \left[ - \ \frac{(x-y+\mu t)^2}{2\sig^2 t} \right] \\
&-& \exp \left[- \frac{2\mu x}{\sig^2} - \frac{(x+y-\mu t)^2}{2\sig^2 t} \right] \; \Bigg\}, \ \ x,y > 0,
\end{eqnarray*}
is the Dirichlet Green's function.  Hence $p(x,t)$ is given by the formula,
\be \label{R2}
p(x,t) = 1 - \int^\infty_0 \ G(x,y,t)dy.
\ee
It is easy to see from this that
\be \label{O2}
\lim_{t \ra 0} \ \sqrt{2\pi \sig^2 t} \ \frac{\pa p}{\pa x}(0,t) = -2 \;,
\ee
whence (\ref{N2}),(\ref{O2}) imply  ${\displaystyle{\lim_{t\ra 0}}} \; \beta(t) = 0$.  One can easily see similarly that ${\displaystyle{\lim_{t\ra \infty}}} \; \beta(t) = \infty$. 
\end{proof}
We may use the function $\beta(t)$ of (\ref{M2}) to improve Proposition 2.1.
\begin{proposition}  Let $\De > 0$ and $\beta(\De) > 0$ satisfy the inequality $K \ge \mu/\rho - u_0 -\beta(\De)$.  Then the function $V(x)$ of (\ref{L2}) satisfies the system of inequalities (\ref{J1}).
\end{proposition}
\begin{proof}  From Lemma 2.2 the function $MV(x)$ of (\ref{P2}) satisfies the inequality $MV(x) \le V_0(x), \; x \ge 0$.  To show that (b) of (\ref{J1}) holds we then need to prove
\[	MV(x) \le x - u_0 + \mu/\rho, \ \ \ x > u_0.	\]
This follows from Lemma 2.1 since $K \ge 0$. 
\end{proof}
Next we consider situations for which $K$ does not satisfy the conditions of Proposition 2.2.  To help us understand this we define a function $u_2(\beta, t)$, $\beta \ge 0$, $t > 0$ as follows:
\[	u_2(\beta, t) = u_0 \ \ \ {\rm if} \ \ \ \beta \le \beta(t).	\]
If $\beta > \beta(t)$ then $u_2(\beta, t)$ is the supremum of all $z \in \R$ such that
\[  \beta\{ 1-p(x,t)\} + h(x,t) \le e^{\rho t} \; V_0(x - z + u_0),  \ x \ge 0.	\]
It is evident that $u_2(\beta, t)$ is a monotonic decreasing function of $\beta$.
\begin{lem} The function $g(\beta) = \beta + u_2(\beta, t), \ \beta > 0$, is strictly monotone increasing.
\end{lem}
\begin{proof}  Suppose $\beta > \beta(t)$.  Since $V'_0(w) \ge 1, \ w \in \R$, we have that for any $\delta > 0$,
\[  e^{\rho t} \; V_0(x - u_2(\beta, t) + \del + u_0) \ge \del e^{\rho t} + \beta\{1-p(x,t)\} + h(x,t), \ x \ge 0.   \]
Hence,
$$ u_2(\beta+\delta e^{\rho t},t)\ge u_2(\beta,t)-\delta,$$
whence $g(\beta)$ is strictly monotonic. 
\end{proof}
\begin{lem} For any $\beta > \beta(t)$ there exists a unique point $u_1(\beta, t) > 0$ such that
\begin{eqnarray*}
\beta\{ 1-p(x,t)\} + h(x,t) &<& e^{\rho t} \; V_0(x - u_2(\beta,t) + u_0), \  x \ge 0, \; x \not= u_1(\beta,t), \\
\beta\{ 1-p(x,t)\} + h(x,t) &=& e^{\rho t} \; V_0(x - u_2(\beta,t) + u_0),  \ x = u_1(\beta,t).
\end{eqnarray*}
\end{lem}
\begin{proof}  Suppose for a given $T > 0$ we have $\beta > \beta(T)$ and define $u(x,t),  \ t>0, x > 0,$ by
\[	u(x,t) = \beta\{ 1-p(x,t)\} + h(x,t) - e^{\rho t} \; V_0(x - u_2(\beta,T) + u_0).  \]
Then $u(x,t)$ is a solution of the diffusion equation (\ref{F2}) with initial and boundary conditions,
\[  u(x,0) = \beta + x - V_0( x-u_2(\beta,T) + u_0), \ u(0,t) = -e^{\rho t}\; V_0(u_0 - u_2(\beta,T)). \]
Evidently $u(x,0)$ is a monotonic decreasing function and ${\displaystyle{\lim_{x\ra \infty}}}\; u(x,0) = -\infty$.  Hence we must have
\[	\beta > V_0( u_0 - u_2(\beta, T) ),	\]
since otherwise $u(x,T) < 0, \; x \ge 0$, which would contradict the definition of $u_2(\beta, T)$.

We show that for small $t$ the set $\{ x>0 : u(x,t) > 0\}$ is an open interval $\big(a(t), b(t) \big)$ with $0 < a(t) < b(t) <\infty$.  To see this first note from (\ref{K2}), (\ref{O2}) that since $p(x,t)$ is a convex decreasing function there is the inequality,
\be \label{Q2}
1 \le \frac{\pa h}{\pa x} \le 1 + 4\mu\sqrt{t}\, / \, \sqrt{2\pi \sig^2}, \ x > 0.
\ee

We can also see from (\ref{R2}) that for any $\del, \; 0 < \del < 1$, there exist positive constants $C(\del)$, $K(\del)$ such that for $0 < t < 1$,
\begin{eqnarray} \label{S2}
p(x,t) &<& \del, \ \ \ x > C(\del)\; \sqrt{t}, \\
\frac{\pa p}{\pa x} (x,t) &<& - \frac{K(\del)}{\sqrt{t}}, \ \ 0 < x < C(\del) \; \sqrt{t}. \nn
\end{eqnarray}
Choose now $\del > 0$ such that
\[  \beta(1 - 3\del) > V_0(u_0 - u_2(\beta,T)).	\]
Then from (\ref{Q2}), (\ref{S2}) there exists $\ve > 0$ such that for $0 < t < \ve$ the function $u(x,t)$ is increasing for $0 < x < C(\del)\; \sqrt{t}$ and $u\big(C(\del)\sqrt{t}, \; t\big) > \beta\del$.

Next we find a region where $u(x,t)$ is decreasing.  To see this observe that we may choose $\ve > 0$ sufficiently small so that for $0 < t < \ve$, there is the inequality,
\[	\frac{\pa p}{\pa x} (x,t) > -\exp\left[ -1/t^{1/6} \right], \ \ x > t^{1/3}.	\]
It follows then from (\ref{K2}) that $u(x,t)$ is decreasing for $x > t^{1/3}$ provided $0 < t < \ve$.  One also has that
\[  \inf \left\{ u(x,t) - u( C(\del)\sqrt{t}, \; t) : C(\del)\sqrt{t} < x < t^{1/3} \right\} \ge - C_1\,t^{1/3},  \]
for some constant $C_1$ which depends on $\del$.  If we choose now $\ve$ small enough so that $C_1 \ve^{1/3} < \beta \del$ then it follows that $\{ x : u(x,t) > 0\} $ consists of an open interval provided $0 < t < \ve$.

We now invoke Theorem A2 of the appendix.  By the definition of $u_2(\beta, T)$ one has $u(x,T) \le 0, \; x \ge 0$, and there is a unique $x = u_1(\beta, T)$ for which $u(x, T) = 0$. 
\end{proof}
\begin{proposition} Let $\De >0$ and $\beta(\De) > 0 $ satisfy the inequality, $K < \mu/\rho - u_0 - \beta(\De)$.  
Then there is a unique solution $\beta > \beta(\De)$ to the equation,
\be \label{T2}
\beta + u_2(\beta, \De) = \mu / \rho - K.
\ee
For this value of $\beta$ put $u_2 = u_2(\beta, \De)$ and  $u_1 = u_1(\beta, \De)$.  Define the function $V(x), x \ge 0$, by
\begin{eqnarray} \label{V2}
V(x) &=& x - u_2 + \mu/\rho, \ \ x > u_2 \ ,  \\
V(x) &=& V_0(x + u_0 - u_2), \  u_1 < x \le u_2 \; , \nn \\
V(x) &=& e^{-\rho \De} \left\{ \beta[1-p(x,\De)] + h(x,\De) \right\}, 0 \le x \le u_1. \nn
\end{eqnarray}
Then the function $V(x)$ satisfies the system of inequalities, (\ref{J1}).
\end{proposition}
\begin{proof}  Since $\beta(\De) + u_2\big( \beta(\De), \De) = \beta(\De) + u_0 < \mu/\rho - K$ it follows from Lemma 2.3 that there is a unique $\beta > \beta(\De)$ satisfying (\ref{T2}).  We can also see that since $K \ge 0$ there is the inequality $u_1 < u_2$.  In fact one has 
\be \label{U2}
\beta [1 - p(x,\De)] + h(x,\De) < e^{\rho \De} V_0(x + u_0 - u_2), \ x \ge u_2.
\ee
This follows from (\ref{K2}) since the LHS of (\ref{U2}) is strictly less than
\[	\mu/\rho -K - u_2 + x + \mu \De < e^{\rho\De} [\mu/\rho + x - u_2] \le e^{\rho \De} V_0(x + u_0 - u_2), \ x \ge u_2,\] 
provided $K \ge 0$.

It is clear now that the function $V$ of (\ref{V2}) is a $C^1$ function and $C^2$ except possibly at the point $x=u_1$.  It is also concave with slope 1 \ for $x \ge u_2$.  Hence (a), (d) of (\ref{J1}) hold.  Next we prove (b).  In view of the concavity of $V$ we have that the function $MV$ of (\ref{K1}) is given by the expression
\be \label{W2}
MV(x) = e^{-\rho \De} \left\{ \left[ \frac \mu \rho - u_2 - K \right] \Big\{1-p(x,\De)\Big\} + h(x,\De)\right\}.
\ee
Thus $V(x) = MV(x), \ 0 \le x \le u_1$, and by the definition of $u_1$ one has $MV(x) \le V(x), u_1 \le x \le u_2$. For $x \ge u_2$ we also have $MV(x) \le V(x)$ by the same argument used to show (\ref{U2}).  We have proved (b) and also (e).

We are left to prove (c).  It is easy to see that $(A-\rho) V(x) \le 0$, $x>u_1$.  We consider then the case $0< x <u_1$.  To see this we observe that the function $u_1(\beta', \De), \ \beta' > \beta(\De)$, is continuous and satisfies ${\displaystyle{\lim_{\beta' \ra \beta(\De)}}} \; u_1(\beta', \De) = 0.$  Hence if $0 < x < u_1$, there exists $\beta'$, $0 < \beta' < \beta$ such that $x = u_1(\beta', \De)$.  Let $w(z)$ be the function,
\[	w(z) = \beta' \{ 1 - p(z,\De)\} + h(z,\De), \ \ z > 0.	\]
By Lemma 2.4 it follows that $(A-\rho)w(z) \le 0$ at $z = u_1(\beta', \De) = x$.  Since
\[ e^{\rho \De} \; MV(z) = (\beta - \beta') \{1 - p(z,\De)\} + w(z),	\]
and $Ap(z,\De) = \pa p/\pa t(z,\De) \ge 0$ one has therefore that $(A-\rho) V(x) \le 0$. 
\end{proof}

\section{Uniqueness of the Solution}

Here we show that the solution to the systems of inequalities (\ref{J1}) is unique provided we make some smoothness assumptions on the function $V(x)$.  Our first goal is to show that a limiting set of strategies realizes the function $V(x)$ constructed in propositions 2.1-2.3. In the following we shall use the convention that if $F(t), \ t>0$, is a right continuous function of time then  at time $\tau$, $F(\tau^+)$ denotes the limit of $F(t)$ as $t$ converges to $\tau$ from above.

We first consider the situation in propositions 2.1-2.2.  Let $\ve$ satisfy $0 < \ve < u_0$.  We define a strategy $\pi_\ve$ for the control process  (\ref{E1}).  Suppose the process begins at $x$ with $0 < x < u_0$.  For those paths which exit the interval $[0,u_0]$ through $u_0$ let $\tau_1$ be the exit time.  We set $L(t)=0, \; t \le \tau_1$, and $L(\tau_1^+) = \ve$.  Thus $X(\tau_1^+) = u_0-\ve$.  If the process begins at $x$ with $x \ge u_0$ we set $\tau_1=0, L(\tau_1^+) = $  $\ve + x - u_0$, whence again $X(\tau_1^+) = u_0-\ve$.  Next we define $\tau_2 > \tau_1$ as the first time the diffusion process with $X(\tau_1^+) = u_0-\ve$ hits $u_0$ for paths which exit the interval $[0, u_0]$ through $u_0$.  We put $L(t) - L(\tau_1^+) = 0$, $\tau_1 < t < \tau_2$, $L(\tau_2 ^+) - L(\tau_1^+) = \ve$.  Thus $X(\tau_2^+) = u_0 - \ve$.  We proceed in this manner defining a sequence of stopping times $\tau_1, \tau_2, \dots$ until the diffusion exits $[0, u_0]$ through 0.
\begin{lem}
Let $V_\ve$ be the return function (\ref{F1}) for the strategy $\pi = \pi_\ve$.  Then $\displaystyle{\lim_{\ve\ra 0}} \; V_\ve(x) = V(x) $ where $V(x)$ is given by (\ref{L2})
\end{lem}
\begin{proof}  Evidently we have that 
\be \label{A3}
V_\ve(x) = x - u_0 + \ve + V_\ve(u_0 - \ve), \ \ x \ge u_0.
\ee
For the diffusion process started at $x, \; 0 < x < u_0$, let $\tau_x$ be the first exit time from the interval $[0,u_0]$.  Then we also have that
\begin{eqnarray} \label{B3}
V_\ve(x) &=& V_\ve(u_0) E\left[ \exp(-\rho \tau_x) ; X(\tau_x) = u_0 \right] \\
&=& V_\ve(u_0) w(x), \ \ 0 < x < u_0 \; , \nn
\end{eqnarray}
where the function $w(x)$ satisfies
\[	(A-\rho) w(x) = 0, \ 0 < x < u_0, \ \ w(0)=0, \ w(u_0) = 1.  \]
It follows that $w(x) = \rho \; V_0(x)/\mu$.  Letting $x = u_0 -\ve$ in (\ref{B3})and using (\ref{A3}) we conclude that $ V_\ve(u_0) $ is given by the formula,
\[	 V_\ve(u_0)  = \ve \; \big/ \; [1 - w(u_0 - \ve)].	\]
Hence $V_\ve(x), \ 0 < x < u_0,$ is given by the formula,
\[     V_\ve(x) = \left[ \frac{ V_0(u_0) -  V_0(u_0-\ve)}{\ve} \right]^{-1} \ V_0(x), \ 0 < x < u_0.	\] 
Since $V'_0(u_0) = 1$ the result follows. 
\end{proof}
Next we consider the situation in Proposition 2.3.  We define a strategy $\pi_\ve$ for the control process (\ref{E1}) whose limiting return function as $\ve \ra 0$ yields the function (\ref{V2}).  If the process begins at $x$ with $x \ge u_2$ we set $\tau_1 = 0, \; L(\tau_1^+) = \ve + x - u_2$, whence $X(\tau_1^+) = u_2 - \ve$.  We require that $0 < \ve < u_2 - u_1$, whence $u_1 < X(\tau_1^+) < u_2$.  If the process begins at $x$ with $u_1 < x < u_2$ we set $\tau_1$ to be the first exit time of the diffusion process from the interval $[u_1,u_2]$.  If $X(\tau_1) = u_2$ then we put $L(t)=0, \; t \le \tau_1$, and $L(\tau_1^+) = \ve$, whence $X(\tau_1^+) = u_2 - \ve$.  Suppose now $X(\tau_1) = u_1$.  We restrict ourselves to all paths of the diffusion process $X(t), \; \tau_1 \le t \le \tau_1 + \De$, which satisfy $X(t) > 0$.  For these paths we set $L(t)=0, \ t \le \tau_1 + \De$,
\be \label{C3}
L((\tau_1 + \De)^+) = X(\tau_1 + \De) - u_2 - K.
\ee
Note that if $\De$ is small the expression in (\ref{C3}) is negative.  We finally put $X((\tau_1  + \De)^+) = u_2.$ For $0 < x \le u_1$ we restrict ourselves to all paths of the diffusion process $X(t), t \le \De$, which satisfy $X(t) > 0$.  For these paths we set $L(t) = 0, \ t \le \De$,
\be \label{D3}
L(\De^+) = X(\De) - u_2 - K.
\ee
Finally we put $X(\De^+) = u_2$.  The process $L(t), t > \tau_1$, is defined similarly.  
\begin{lem}
Let $V_\ve$ be the return function (\ref{F1}) for the strategy $\pi = \pi_\ve$.  Then $\displaystyle{\lim_{\ve\ra 0}} \; V_\ve(x) = V(x) $ where $V(x)$ is given by (\ref{V2}).
\end{lem}
\begin{proof}  Arguing as in Lemma 3.1 we have that
\be \label{E3}
V_\ve(x) = x - u_2 + \ve + V_\ve(u_2 - \ve), \ \ x \ge u_2.
\ee
\be \label{F3}
V_\ve(x) =V_\ve(u_2) w_2(x) + V_\ve(u_1)w_1(x), \ \ u_1 < x < u_2,
\ee
where
\begin{eqnarray*}
w_2(x) &=& E\left[ \exp(-\rho \tau_x); \ X(\tau_x) = u_2 \right], \\
w_1(x) &=& E\left[ \exp(-\rho \tau_x); \ X(\tau_x) = u_1 \right],
\end{eqnarray*}
and $\tau_x$ is the exit time from the interval $[u_1, u_2]$ for the diffusion process started at $x$.  For $0 < x \le u_1$ we have in addition the identity,
\be \label{G3}
V_\ve(x) = e^{-\rho \De} E\left[ X(\De) - u_2 - K \; ; \; \tau_x > \De \right] + e^{-\rho \De} V_\ve(u_2) P(\tau_x > \De),
\ee
where $\tau_x$ is the first time the diffusion started at $x$ hits 0.  From (\ref{V2}) we can rewrite (\ref{G3}) as 
\be \label{H3}
V_\ve(x) = V(x) + e^{-\rho \De} \big\{ V_\ve(u_2) - V(u_2)\big\} [1 - p(x , \De)].
\ee
It is clear that for $u_1 < x < u_2$ the function $V(x)$ may be written as
\[   V(x) = V(u_2) w_2(x) + V(u_1)w_1(x), \ \ u_1 < x < u_2.	\]
Hence if we put $g_\ve(x) = V_\ve(x) - V(x)$ we have from (\ref{F3}) that
\be \label{I3}
g_\ve(x) =g_\ve(u_2) w_2(x) + g_\ve(u_1)w_1(x), \ \ u_1 < x < u_2.
\ee
Setting $x=u_1$ in (\ref{H3}) we also have that 
\be \label{J3}
g_\ve(u_1) =  e^{-\rho \De} g_\ve(u_2) [1 - p(u_1 , \De)].
\ee
Putting $x = u_2$ in (\ref{E3}) yields the identity,  
\be \label{K3}
g_\ve(u_2 - \ve) = g_\ve(u_2) + [V(u_2) - V(u_2 - \ve) - \ve].
\ee
If we set $x=u_2 - \ve$ in (\ref{I3}) and use (\ref{J3}), (\ref{K3}) we obtain a formula for $g_\ve(u_2)$,
\be \label{L3}
g_\ve(u_2) = - [V(u_2) - V(u_2 - \ve) - \ve] \Big/ \left\{ 1 - w_2(u_2-\ve) - e^{-\rho \De}  [1 - p(u_1 , \De)] w_1(u_2-\ve)\right\}.
\ee
Since $V$ is concave at $u_2$ the numerator of (\ref{L3}) is negative.  The denominator is positive since $w_1(x) + w_2(x) < 1$, $u_1 < x < u_2$.  Hence $g_\ve(u_2) < 0$.  It follows now from (\ref{E3}), (\ref{H3}), (\ref{I3}) that $V_\ve(x) < V(x), \ x > 0$.  Since the numerator of (\ref{L3}) is $O(\ve^2)$ and the denominator  is from Hopf's maximum principle bounded below by a positive constant times $\ve$, it follows that  ${\displaystyle \lim_{\ve \ra 0}} \  g_\ve(u_2) = 0$.  Hence  ${\displaystyle \lim_{\ve \ra 0}} \  V_\ve(x) = V(x), \ x \ge 0$. 
\end{proof}
We have shown that certain limiting strategies yield the return functions given in propositions 2.1 -2.3.  Next let $V(x)$ be a $C^1$ solution of the system of inequalities (\ref{J1}).  Since $V'(x) \ge 1, \ x > 0$, the limit,
\be \label{M3}
\beta + K = \lim_{x \ra \infty} \ [V(x) - x] \ \ {\rm exists}.
\ee
This limit must be finite.  Otherwise the function $MV(x)$ cannot be finite.  Hence from (\ref{K1}) we have
\be \label{N3}
MV(x) = e^{-\rho \De} \big\{ \beta [1 - p(x , \De)] + h(x, \De)\big\}.
\ee
\begin{lem}  There exists $u_2, \ve$ with $0 < \ve < u_2$ such that $V(x) = \beta + K + x$ for $x \ge u_2$, and for $u_2 - \ve < x \le u_2$, $V(x)$ is the solution to the initial value problem,
\be \label{O3}
(A - \rho) V(x) = 0, \ \ V(u_2) = \mu /\rho, \ \ V'(u_2) = 1.
\ee
Further, $\beta$ and $u_2$ are related by the identity,
\be \label{P3}
\beta = \mu/\rho - K - u_2.
\ee
\end{lem}
\begin{proof}  Suppose $u_2 > 0$ is a point which has the property that for some $\del > 0$ one has $V'(x)=1$ for $u_2 \le x < u_2 + \del$, and $(A-\rho)V(x) = 0$ for $u_2 - \del < x < u_2$.   Since $(A-\rho)V(x) \le 0$ for $u_2 \le x < u_2 + \del$ it follows that $V(u_2) \ge \mu /\rho$.  Using the fact that $V$ is $C^1$ at $u_2$ and $(A-\rho)V(x) = 0$, $x < u_2$, we can conclude now that $\displaystyle{\lim_{x\ra u_2-}} V''(x) \ge 0$.  If $\displaystyle{\lim_{x\ra u_2-}}V''(x) > 0$ then $V'(x) < 1$ for $x < u_2$ with $u_2 - x$ sufficiently small, in contradiction to (\ref{J1}).  Hence $\displaystyle{\lim_{x\ra u_2-}} V''(x) = 0$, whence $V(u_2) = \mu/\rho$ since $V'(u_2) = 1$.  Thus for $u_2-\delta<x<u_2$ the function $V$ is the solution to (\ref{O3}).

From (\ref{K2}), (\ref{M3}), (\ref{N3}) we see that there exists $u_3 > 0$ such that $MV(x) < V(x)$ for $x \ge u_3$.  Hence for each $x > u_3$ the function $V$ must satisfy $(A-\rho)V(x) = 0$ or $V'(x) = 1$.  Observe now that (\ref{M3}) implies that $\{ x : (A-\rho)V(x) = 0\}$ does not include a neighborhood of $\infty$.  Hence by the argument of the previous paragraph $\{ x : V'(x) = 1\}$ does include a neighborhood of $\infty$.  We define $u_2$ by
\[	u_2 = \inf \{ z : V'(x) = 1, \ x \ge z \}.	\]
We show that $MV(u_2) < V(u_2)$.  To see this first note that from (\ref{M3}) one has $V(x) = \beta +K+x, \ x > u_2$, and we also have that $V(u_2) \ge \mu/\rho$.  Now $V'(x) \ge 1, \; 0 < x < u_2$, and $V(0) = 0$.  Hence $\beta + K \ge 0$.  If $\beta \le 0$ then (\ref{K2}), (\ref{N3}) yield the inequality,
\[	MV(u_2) < e^{-\rho \De} (u_2 + \mu \De).	\]
On the other hand we have
 \[	e^{\rho \De} V(u_2) \ge \left[ e^{\rho \De}-1\right] \mu/\rho+ \beta + K + u_2 \ge \mu \De + u_2 \ ,	\] 
whence $MV(u_2) < V(u_2).$  If on the other hand $\beta \ge 0$, then
\begin{eqnarray*}
	MV(u_2) &<& e^{-\rho\De} [\beta + u_2 + \mu \De], \\
e^{\rho\De} V(u_2) &\ge& \mu\De + \beta + K + u_2 .
\end{eqnarray*}
Hence again we have $MV(u_2) < V(u_2)$.

The result of the lemma now easily follows since by the previous paragraph there exists $\ve > 0$ such that $(A-\rho)V(x) = 0$ for $u_2 - \ve < x < u_2$.  From the first paragraph it follows that $V$ is the solution to (\ref{O3}).  The identity (\ref{P3}) follows from the fact that $V(u_2) = \mu /\rho$. 
\end{proof}
Next we define $u_1 < u_2$ by
\[	u_1 = \inf \{ z > 0 : (A - \rho) V(x)=0, \ z < x <u_2\}.	\]
\begin{lem}  If $u_1 > 0$ then $u_2 < u_0, \; \beta > 0$ and $V(x) = MV(x)$ for $0 < x < u_1$.
\end{lem}
\begin{proof}  We proceed as in Lemma 2.4 by considering the function $u(x,t)$ given by 
\[	u(x,t) = \beta \{ 1 - p(x,t)\} + h(x,t) - e^{\rho t} V_0(x - u_2 + u_0).	\]
Then $u(x,t)$ is a solution of the diffusion equation (\ref{F2}) with initial and boundary conditions,
\[	u(x,0) = \beta + x - V_0(x-u_2 + u_0), \ u(0,t) = - e^{\rho t} V_0(u_0 - u_2 ).\]
Let us suppose first that  $u_2 \ge u_0$.  In that case $u(0,t) \ge 0$ and $u(x,0)$ is a monotonic decreasing function with $\displaystyle{\lim_{x\ra \infty}} u(x,0) = -\infty$.  It is easy to see from this that for small $t$ the function $u(x,t)$ has at most one sign change.  Hence by Theorem A1, $u_1$ is the unique solution to the equation $u(x,\De)=0$. From Theorem A3 it follows that $\pa u/\pa x (x, \De) < 0$ at $x = u_1$, but this contradicts the $C^1$ property of the function $V(x)$ at $x=u_1$.  We conclude that $u_2 < u_0$.

Assuming $u_2 < u_0$, then $u(0,t) < 0$.  Hence $\beta > V_0(u_0 - u_2) > 0$ since otherwise $u(x,t) < 0, \; x \ge 0, \; t>0$.  By Theorem A2 the set $\{ x > 0 : u(x, \De) \ge 0\}$ is a closed interval with $u_1$ as one of its end points.  Evidently $u_1$ must be the right most end point.  If the interior of the interval is non empty then $\pa u/\pa x (x, \De) < 0$ at $x = u_1$ by Theorem A3 of the appendix.  Since this again contradicts the $C^1$ property of  $V$ at $u_1$ we conclude that $\{ x > 0 : u(x,\De) \ge 0\} = \{u_1\}$.  It is clear now that in the notation  of $\S 2$ we have $\beta > \beta(\De)$, $u_2 = u_2(\beta, \De)$ and $u_1 = u_1(\beta, \De)$.

Finally we need to show that $V(x) = MV(x), \; 0 < x < u_1$.  Let $u_3 = \inf\{ z : 0 < z < u_1$,  $V(x) = MV(x)$ for $z < x < u_1\}$.  Evidently $0 \le u_3 < u_1$.  Suppose now $u_3 > 0$.  Since $V(x)$ is concave for $x > u_3$ it follows that $V'(u_3) > 1, \; V(u_3) < \mu / \rho$.  Let $u_4 = \inf\{ z : 0 < z < u_3$,  $(A-\rho)V(x) = 0$ for $z < x < u_3\}$.  It is easy to see that $V(x)$ is concave for $u_4 < x < u_2$.  Let $w(x) = V(x) - MV(x), \  u _4 < x < u_3$.  Then we must have that $w(u_4) = w(u_3) = 0$ and $w(x) \ge 0,  \ u_4 < x < u_3$.  By the argument of Proposition 2.3 we also have that $(A-\rho)w(x) \ge 0, \  u_4 < x < u_3$.  Hence by the maximum principle it follows that $w(x) = 0,  \ u_4 < x < u_3$.  Since this contradicts the definition of $u_3$ we conclude that $u_3 = 0$. 
\end{proof}
\begin{proposition}  Let $V(x),  \ x \ge 0$ be a $C^1$ solution to the set of inequalities (\ref{J1}).  Then $V$ is the unique solution given by Propositions 2.1-2.3.
\end{proposition}
\begin{proof}  This follows from Lemmas 3.3, 3.4. 
\end{proof}
We give an alternative proof of Proposition 3.1 which avoids the use of Theorem A3.  Instead, we shall use the technique of ``verification theorem", \cite{fr,fls}
\begin{lem}  With $u_1, u_2$ as defined in Lemma 3.4 there is the inequality $u_2 \ge u_2(\beta, \De)$, where $\beta$ is the solution to (\ref{T2}).
\end{lem}
\begin{proof}  Let $\pi_\ve$ be the strategy defined just before Lemma 3.2, where $u_1, u_2$ are as in Lemma 3.4.  If $V_\ve$ is the return function (\ref{F1}) ) corresponding to $\pi_\ve$, then by the argument of Lemma 3.2 we see that ${\displaystyle \lim_{\ve \ra 0}} \  V_\ve(x) = V(x)$ for $x \ge u_1$, where $V(x)$ is the function discussed in Lemma 3.3 and Lemma 3.4.

Next let $V_{opt}(x)$ be the solution to the control problem constructed in propositions 2.1-2.3.  We shall show that
\be \label{Q3}
V_{opt}(x) \ge V_\ve(x), \ \ x \ge u_1.
\ee
To see this first let $X(t),  \ t > 0$, be the diffusion process with generator $A$.  Assume $X(0) > 0$ and $\tau$ is the first hitting time at 0.  Since $V_{opt}(x)$ is $C^1$ for $x \ge 0$ and $C^2$ for all $x \ge 0$ except possibly $x = u_1(\beta, \De)$ with $\beta$ satisfying (\ref{T2}), it follows that
\be \label{R3}
e^{-\rho t} V_{opt}\left( X(t \wedge\tau)\right) - \int^{t\wedge\tau}_0 e^{-\rho t} (A-\rho)V_{opt}\left( X(s)\right) ds
\ee
is a martingale.

Now let $X(0) = u_2 - \ve$ and $\tau_1$ be the first exit time from the interval $[u_1, u_2]$.  Evidently $\tau_1 < \tau$.  Since $(A-\rho)V_{opt}(x) \le 0, \; x \ge 0$ it follows from (\ref{R3}) that 
\be \label{S3}
V_{opt}(u_2 - \ve) \ge E\left[ e^{-\rho \tau_1} V_{opt}\left( X(\tau_1)\right) \right].
\ee
Since $V'_{opt}(x) \ge 1, \ x \ge 0$, it follows that 
\be \label{T3}
E\left[ e^{-\rho \tau_1} V_{opt}\left( X(\tau_1)\right); X(\tau_1) = u_2 \right] \ge
\ee
\[ E\left[ e^{-\rho \tau_1} V_{opt}(u_2 -\ve) ; X(\tau_1) = u_2 \right]   
 \  + E\left[ e^{-\rho \tau_1} \ve; X(\tau_1) = u_2 \right] .	\]
On using the fact that $MV_{opt}(x) \le V_{opt}(x), \; x \ge 0$, we also have that
\begin{eqnarray*}
&\  & E\left[ e^{-\rho \tau_1} V_{opt} \left( X(\tau_1)\right);  X(\tau_1) = u_1 \right] \\
&\ge& E \left[ e^{-\rho( \tau_1+\De)} V_{opt}(u_2)  ;  X(\tau_1) = u_1, \; \tau_1 + \De < \tau \right] \\
&+& E\left[ e^{-\rho( \tau_1+\De)} \left[ X(\tau_1 + \De) - u_2 - K\right]  ; X(\tau_1) = u_1, \tau_1 + \De < \tau\right]. 
\end{eqnarray*}
Hence there is the inequality,
\begin{eqnarray} \label{U3}
& \ & E\left[ e^{-\rho \tau_1} V_{opt} \left( X(\tau_1)\right);  X(\tau_1) = u_1 \right] \\
&\ge& E\left[ e^{-\rho( \tau_1+\De)} V_{opt}(u_2-\ve)  ;  X(\tau_1) = u_1, \; \tau_1 + \De < \tau \right] \nn \\
&+& E\left[ e^{-\rho( \tau_1+\De)} \ve ; X(\tau_1) = u_1, \; \tau_1 + \De < \tau \right] \nn \\
&+& E\left[ e^{-\rho( \tau_1+\De)} \left[ X(\tau_1 + \De) - u_2 - K\right]  ; X(\tau_1) = u_1, \tau_1 + \De < \tau\right] \nn.
\end{eqnarray}
If we define now $\tau^*_1$ as $\tau^*_1 = \tau_1$ if $X(\tau_1) = u_2, \; \tau^*_1=\tau_1 + \De$ if $X(\tau_1) = u_1$ we have from (\ref{S3}), (\ref{T3}), (\ref{U3}), the inequality,
\be \label{V3}
V_{opt}(u_2 - \ve) \ge E\left[ e^{-\rho \tau_1^*} V_{opt}(u_2-\ve)  ;   \tau_1^* < \tau \right] 
\ee
\[ + E\left[ \int^{\tau^*_1}_0 \ e^{-\rho t} dL(t) ; \tau^*_1 < \tau \right],\]
where $L(t)$ is the return function associated with the strategy $\pi_\ve$.  Evidently if we iterate the inequality (\ref{V3}) we obtain (\ref{Q3}) for $x = u_2 - \ve$.  The inequality for all $x \ge u_1$ follows in a similar way.

If we let $\ve \ra 0$ in (\ref{Q3}) we obtain the inequality $V_{opt}(x) \ge V(x), \; x \ge u_1$, which implies the result. 
\end{proof}
\begin{lem} With $u_1,u_2$ as defined in Lemma 3.4 there is the inequality $u_2 \le u_2(\beta, \De)$, where $\beta$ is the solution to (\ref{T2}).
\end{lem}
\begin{proof} Let $\pi_{opt, \ve}$ be the strategy of Lemma 3.2 and $V_{opt, \ve}$ the corresponding return function.  Then by Lemma 3.2 we have that ${\displaystyle \lim_{\ve \ra 0}}  \ V_{opt,\ve}(x) = V_{opt}(x)$ where $V_{opt}$ is the function given by (\ref{V2}).  We shall show that for $V$, the function discussed in Lemma 3.3 and Lemma 3.4, there is the inequality,
\be \label{W3}
V(x) \ge V_{opt, \ve}(x), \ \ \ x \ge u_1(\beta, \De),
\ee
where $\beta$ is the solution to (\ref{T2}).  In fact the proof of (\ref{W3}) is identical to the proof of (\ref{Q3}) since $V$ satisfies the variational inequalities (\ref{J1}).  The result follows by letting $\ve \ra 0$. 
\end{proof}
\begin{proof}[Proof of Proposition 3.1]  We have shown in Lemma 3.5 and Lemma 3.6 that $u_2 = u_2(\beta, \De)$ with $\beta$ satisfying (\ref{T2}).  Since $MV$ is given by (\ref{N3}) and $MV(u_1) = V(u_1)$ we must have $u_1 = u_1(\beta,\De)$.  The fact that $V(x) = MV(x)$ for $0 < x < u_1$ follows by the argument at the end of Lemma 3.4. 
\end{proof}

\section{Properties of the Thresholds $u_1, \ u_2$.}
In this section we shall study the properties of $u_1,u_2$ as defined in Proposition 2.3.  Evidently $u_1,u_2$ are functions of $K \ge 0$ and $\De > 0$.  If $K \ge \mu/\rho - u_0$ then $u_2 = u_0, u_1 = 0$.  We shall therefore be interested in the situation where $0 \le K < \mu/\rho - u_0$.
\begin{lem}  Suppose $0 \le K < \mu/\rho - u_0$ and for $\De > 0$ let $u_{1,K}(\De)$, $u_{2,K}(\De)$ be the values  of  $u_1,u_2$ determined by $K,\De$.  Then $u_{1,K}, u_{2,K}$ are continuous functions satisfying:
\be \label{A4}
\lim_{\De \ra 0} u_{1,K}(\De) = 0, \ \ \ \lim_{\De \ra 0} u_{2,K}(\De) = u_0 - \hat u_K \ ,
\ee
where $z = \hat u_K$ is the unique solution to the equation,
\be \label{B4}
V_0(z) = z + \mu/\rho - K - u_0.
\ee
$u_{1,K}(\De) = 0$, $u_{2,K}(\De) = u_0$, provided $\De \ge \De_0$ where $\De = \De_0$ is the unique solution to the equation
\[	\beta(\De) = \mu/\rho - K - u_0.	\]
\end{lem}
\begin{proof}  It is easy to see that $u_{1,K}(\De)$ and $u_{2,K}(\De), \ \De > 0$, are continuous functions.  Evidently Proposition 2.2 implies that $u_{1,K}(\De) = 0$, $u_{2,K}(\De) = u_0$ if $\De \ge \De_0$.  We consider the case $\De \ra 0$.  Then the function $u_2(\beta, \De)$ defined just before Lemma 2.3 satisfies 
\[	\lim_{\De \ra 0} u_2(\beta, \De) = u_0 - z_\beta\ ,	\]
where $z_\beta$ is the unique solution to the equation $V_0(z_\beta) = \beta$.  From Proposition 2.3 it follows therefore that $\hat u_K = z_\beta$ where $\beta$ satisfies $\beta + u_0 - z_\beta = \mu/\rho - K$.  This equation is evidently the same as $(\ref{B4})$. Note that $0 < \hat u_K \le u_0$ since $V_0(0) = 0, \; V_0(u_0) = \mu/\rho$ and $V'_0(z) > 1, \; z > 0$. 
\end{proof}
Next we obtain the first order behavior of $\beta(\De), u_{1,K}(\De),u_{2,K}(\De)$ as $\De \ra 0$.  The first order behavior of $\beta(\De)$ can easily be obtained from (\ref{M2}), (\ref{O2}), (\ref{Q2}).  Thus we have 
\[	\lim_{\De \ra 0} \beta(\De) / \sqrt{\De} = \left[ V'_0(0) - 1\right] \sqrt{\pi\sig^2/2} >0.	\]
To obtain the first order behavior of $u_{1,K}( \De), u_{2,k} (\De)$ as $\De \ra 0$ we consider the behavior of the functions $u_1(\beta, \De), u_2(\beta, \De)$ as $\De \ra 0$ where $\beta > 0$.
\begin{lem}  Let $u_1(\beta, \De), u_2(\beta, \De)$ be the functions defined in Lemmas 2.3, 2.4.  Then if $\beta > 0$ there are the limits ,
\begin{eqnarray*}
\displaystyle{\lim_{\De \ra 0}}  \  u_1 (\beta,\De) / \sig \sqrt{\De} | \ell n \; \De|^{1/2} &= &1, \\
\displaystyle{\lim_{\De \ra 0}}  \  \left\{ u_2 (\beta, \De) - [u_0 - \hat u_\beta]\right\} \Big/ \sig \sqrt{\De} |
 \ell n \; \De|^{1/2} &=& 1 - 
1 / V'_0(\hat u_\beta),
\end{eqnarray*}
where $\hat u_\beta$ is the unique solution $z= \hat u_\beta$ to the equation $V_0(z) = \beta$. 
\end{lem}
\begin{proof}  Just as in the proof of Lemma 4.1 we see that $u_1(\beta, \De), u_2(\beta, \De)$satisfy
\be \label{E4}
\lim_{\De \ra 0} u_1(\beta,\De) = 0, \ \ \ \lim_{\De \ra 0} u_2(\beta, \De) = u_0 - \hat u_\beta \ .
\ee
For $\De > 0$ it follows from Lemma 2.4 that $u_1=u_1(\beta, \De)$ and $u_2=u_2(\beta, \De)$ are the unique positive solutions to the system of equations,
\be \label{C4}
\beta \{ 1-p(u_1, \De)\} + h(u_1,\De) = e^{\rho\De} \; V_0(u_1 - u_2 + u_0),
\ee
\be \label{D4}
- \beta \; \frac{\pa p}{\pa x} (u_1, \De) + \frac{\pa h}{\pa x} (u_1, \De) = e^{\rho\De} \; V'_0(u_1 - u_2 + u_0).
\ee
We shall look for solutions to (\ref{C4}), (\ref{D4}) which satisfy (\ref{E4}).  To do this we first note from (\ref{R2}) that $\pa p/\pa x$ is given by the formula,
\be \label{F4}
-\frac{\pa p}{\pa x} \; (x,t) = \frac 2{\sqrt{2\pi\sig^2t}} \ \bigg\{ \exp \left[ - \frac{(x+\mu t)^2}{2\sig^2 t}\right]+
\ee
\[ \frac\mu{\sig^2} \exp\left[ - \frac{2\mu x}{\sig^2}\right] \int^\infty_0 \exp \left[ - \frac{(x+y-\mu t)^2} {2\sig^2 t}\right]dy\bigg\}.	\]
Let $g(z), z > 0$, be the function,
\be \label{L4}
g(z) = \sqrt{2\pi\sig^2\De} \ \left\{ - \beta\; \frac{\pa p}{\pa x} \left( z\sqrt{\De}, \; \De\right) + \frac{\pa h}{\pa x} \left( z\sqrt{\De}, \; \De\right) - e^{\rho \De} \; V'_0 \left( z\sqrt{\De} - u_2 + u_0\right) \right\}.
\ee
where $u_2$ is a fixed  parameter restricted to lie in the region,
\be \label{H4}
\hat u_{\beta}/2 < u_0 - u_2 < ( u_0 + \hat u_\beta)/2.
\ee
Observe now that in view of (\ref{H4}) there are constants $\De_0, K_0 > 0$ such that if $0 < \De < \De_0$ and $a(\De), b(\De) > 0$ are defined by the identities,
\begin{eqnarray} \label{M4}
a(\De)^2 / 2\sig^2 &=& - \frac 1 2 \ln \De - K_0 \; , \\
b(\De)^2 / 2\sig^2 &=& - \frac 1 2 \ln \De + K_0 \; , \nn
\end{eqnarray}
then the function $g$ is strictly monotonic decreasing in the interval $[a(\De), b(\De)]$ with $g(a(\De)) > 0, \; g(b(\De)) < 0$.  It follows that there is a unique solution $z = z_\De(u_2)$ of the equation $g(z) = 0$ in the interval $(a(\De), b(\De))$.  We have shown then that (\ref{D4}) gives $u_1 = \sqrt{\De} \; z_\De(u_2)$ as a unique function of $u_2$ provided $0 < \De < \De_0$ and $u_2$ satisfies(\ref{H4}).

Next we wish to estimate the LHS of (\ref{C4}) when $u_1 = \sqrt{\De}\; z$ and $z \in [a(\De), b(\De)]$.  To do this we write
\[	p\big( \sqrt{\De}\; z, \De\big) =- \int^\infty_{\sqrt{\De}\; z} \frac{\pa p}{\pa x} \; (x,\De)dx,  \]
and use the formula (\ref{F4}).  Observe now that
\[	\frac 2{\sqrt{2\pi\sig^2\De}} \int^\infty_{\sqrt{\De}\;z} \exp \left[ - \frac{(x+\mu\De)^2}{2\sig^2\De} \right]dx \]
\[  = \exp \left[ - \frac 1{2\sig^2} \left( z + \mu\sqrt{\De} \right)^2 \right] \frac 2{\sqrt{2\pi\sig^2}} \int^\infty_0 d\xi  \  \exp \left[-\frac{\xi(z+\mu\sqrt{\De})}{\sig^2} - \frac{\xi^2}{2\sig^2}\right]		\]
\[ \le \exp \left[ - \frac 1{2\sig^2} \left( z + \mu\sqrt{\De} \right)^2 \right] \frac {2\sig}{\sqrt{2\pi}(z+\mu\sqrt{\De})}.\]
We can similarly estimate the contribution to $p(\sqrt{\De}\;z, \De)$ from the second term in (\ref{F4}).  Thus we have
\[	\frac 2{\sqrt{2\pi\sig^2\De}} \int^\infty_0 \exp \left[ - \frac{(x+y-\mu \De)^2} {2\sig^2\De}\right] dy \]
\[   \le \frac{2\sig\sqrt{\De}}{\sqrt{2\pi} (x-\mu\De)}\ \exp \left[- \frac{1}{2\sig^2\De} (x-\mu\De)^2 \right] ,  \]
\[   \int^\infty_{\sqrt{\De}\; z} \; dx  \  \frac{2\sig\sqrt{\De}}{\sqrt{2\pi} (x-\mu\De)} \frac{\mu}{\sig^2} \exp \left[ -\frac{2\mu x}{\sig^2}\right] \exp \left[- \frac{1}{2\sig^2\De} (x-\mu\De)^2 \right] \]
\[   \le \frac{2\sig\sqrt{\De}}{\sqrt{2\pi} (z + \mu\sqrt{\De})}\exp \left[- \frac{1}{2\sig^2} (z+\mu\sqrt{\De})^2 \right] ,  \]
provided $z \ge 2\mu\max(1,\sqrt{\De_0})$.  We conclude then that for $\De_0$ sufficiently small there is a constant $C$ such that 
\be \label{G4}
0 < p\Big( \sqrt{\De}\; z, \De\Big) \le C\sqrt{\De} / |\ln \; \De|^{1/2},
\ee
provided $0 < \De \le \De_0, \ z \in [a(\De), b(\De)]$, and $u_2$ satisfies (\ref{H4}).

Consider now the function $F(u_2)$ defined by 
\[	F(u_2) = e^{\rho \De} \; V_0(u_1 - u_2 + u_0) - h(u_1, \De) - \beta \{ 1-p(u_1,\De)\},  \ 0 < \De < \De_0\; ,  \]
for $u_2$ in the region (\ref{H4}) and $u_1 = \sqrt{\De} \; z_\De(u_2)$ the unique solution of (\ref{D4}), $z_\De(u_2) \in (a(\De), b(\De))$.  In view of (\ref{G4}) it is clear that
\[	F(u_0 - \hat u_\beta) > 0, \ F\left( u_0 - \hat u_\beta + \sqrt{\De} \; b(\De) \right) < 0,		\]
whence there is a solution $u_2$ to the equation $F(u_2) = 0$ in the region $u_0 - \hat u_\beta < u_2 < u_0 - \hat u_\beta + \sqrt{\De} \; b(\De)$.  We have therefore shown the existence of a solution $(u_1, u_2)$ to the set of equations (\ref{C4}), (\ref{D4}) provided $0 < \De < \De_0$.  By Lemma 2.4 the solution is unique.  One can easily now derive the asymptotics of $u_1(\beta, \De), u_2(\beta, \De)$ as $\De \ra 0$.  In fact the asymptotics of $u_1(\beta,\De)$ is already a consequence of the fact that $z_\De(u_2) \in (a(\De), b(\De))$.  To obtain the asymptotics of $u_2(\beta, \De)$ we do a Taylor expansion of $F(u_2)$ about $u_2 = u_0 - \hat u_\beta$.  Thus we have
\[ F(u_2) = e^{\rho \De} \left\{ V_0(\hat u_\beta) +[u_1 - u_2 + u_0- \hat u_\beta] V'_0(\hat u_\beta) + O\left([u_1 - u_2 + u_0- \hat u_\beta]^2\right) \right\} \]
\[	-u_1 - \beta + O \left( \sqrt{\De} / |\ln \; \De|^{1/2}\right).	\]
Hence the asymptotic form of $u_2 = u_2(\beta,\De)$ is obtained from the equation
\[          [u_1 - u_2 + u_0- \hat u_\beta] V'_0(\hat u_\beta) - u_1 = 0,   \]
where $u_1 = u_1(\beta,\De)$ has the asymptotic form $u_1(\beta,\De) = \sig\sqrt{\De} / |\ln \; \De|^{1/2}$. \end{proof}
\begin{proposition}
Let $u_{1,K}(\De), \; u_{2,K}(\De)$ be the functions of $\De$ defined in Lemma 4.1.  Then there are the limits,
\be \label{I4}  
\lim_{\De \ra 0} u_{1,K}(\De) / \sig\sqrt{\De}  \ |\ln \; \De|^{1/2} = 1,
\ee
\be \label{J4}  
\lim_{\De \ra 0} \left\{ u_{2,K}(\De) - [u_0 - \hat u_K] \right\} / \sig \sqrt{\De}  \ |\ln \; \De |^{1/2} = 1,
\ee
where $\hat u_K$ is as in Lemma 4.1.
\end{proposition}
\begin{proof} We use Lemma 4.2 and (\ref{T2}).  Evidently (\ref{I4}) follows directly from Lemma 4.2.  To get (\ref{J4}) we substitute from Lemma 4.2 the formula for $u_2(\beta, \De)$ and solve approximately for $\hat u_\beta$.  Thus on writing $\hat u_\beta = \hat u_K + \del$ we have to highest order,
\begin{eqnarray*}
V_0(\hat u_K + \del) &+& u_0 - \hat u_K - \del \\
&+& [1-1/V'_0(\hat u_K)] \sig\sqrt{\De}  \ |\ln \; \De|^{1/2} = \mu/\rho - K.
\end{eqnarray*}
Taylor expanding this last identity about $\del = 0$ and solving for $\del$ yields,
\[	\del = - \sig\sqrt{\De}  \  |\ln \; \De|^{1/2} / V'_0(\hat u_K).		\]

Hence we have to highest order,
\begin{eqnarray*}
u_{2,K} (\De) &=& [u_0 - \hat u_ \beta] + \sig \sqrt{\De} \  |\ln \; \De |^{1/2} \{ 1 - 1 / V'_0(\hat u_\beta) \} \\
&=& [u_0 - \hat u_K - \del ] + \sig \sqrt{\De}  \  |\ln \; \De|^{1/2} \{ 1 - 1 / V'_0(\hat u_K) \} \\
&=& u_0 - \hat u_K  + \sig \sqrt{\De}  \  |\ln \; \De|^{1/2}. 
\end{eqnarray*}
\end{proof}
Finally we wish to show that the function $V(x), x \ge 0$, of Proposition 2.3, which is a $C^1$ function, fails to be twice differentiable at $x = u_1$.  To do this let $u_2$ satisfy $0 < u_2 < u_0$ and $\beta > V_0(u_0 - u_2)$.  We consider the function $u(x,t)$ defined by
\be \label{K4}
u(x,t) = \beta \{ 1 - p(x,t) \} + h(x,t) - e^{\rho t} V_0(x - u_2 + u_0).
\ee
Evidently $u(x,0) = \beta + x - V_0(x - u_2 + u_0)$ and $u(0,t) = -e^{\rho t} \; V_0(u_0 - u_2),$ whence $u(0,t) < 0, \; t > 0$.
\begin{lem} Let $u(x,t)$ be the function (\ref{K4}).  Then there exists $\ve > 0$ and $x(t) > 0,  \ 0 < t < \ve$, such that $\pa u/\pa x(x,t) > 0$ for $0 \le x < x(t)$ and $\pa u/\pa x(x,t) < 0$ for $x > x(t)$.
\end{lem}
\begin{proof}  We proceed as in Lemma 4.2 observing that the function $g$ of (\ref{L4}) is given by
\be \label{N4}
g(z) = \sqrt{2\pi \sig^2\De} \ \pa u/\pa x \Big( z \sqrt{\De}, \ \De\Big).
\ee
With $a(\De), b(\De)$ defined by $(\ref{M4})$ we have seen that provided $0 < \De < \De_0$ then $g(z) > 0$ for $0 < z \le a(\De)$, $g(z) < 0$ for $z \ge b(\De)$ and $g$ is strictly monotonic decreasing in the interval $[a(\De), b(\De)]$. 
\end{proof}
\begin{lem} Let $u(x,t)$ be the function (\ref{K4}).  Then there exists $\ve > 0$ and $a(t), b(t) > 0,  \ 0 < t < \ve$, such that $\{x > 0 : \pa^2 u/ \pa x^2(x,t) > 0\}$ $= (a(t), b(t))$.
\end{lem}
\begin{proof}  From (\ref{K2}), (\ref{K4}) we have that
\[	\pa^2 u/\pa x^2(x,t) = -\beta \pa^2 p/\pa x^2(x,t) - \mu \int^t_0 \pa^2 p/\pa x^2(x,s)ds - e^{\rho t} \; V''_0(x-u_2 + u_0).	\]
Since $p(x,t)$ is a convex function of $x$ it follows that $\pa^2 u/\pa x^2(x,t) < 0$ for $x > u_2$.  Since $V''_0(u_0) = 0$, $V'''_0(u_0) = 2\rho/\sig^2$ we conclude that there exists $\De_0 > 0$ and $b(t)>0, \; 0 < t < \De_0$, with the property that $\displaystyle{\lim_{t\ra 0}} \; b(t) = u_2$, and $\pa^2 u/\pa x^2(x,t) < 0$ if $x > b(t)$ and $\pa^2 u/\pa x^2(x,t) > 0$ if $t^{1/3} < x < b(t)$.

We consider now the interval $0 < x < t^{1/3}$.  From (\ref{N4}) we may consider the function $g'(z)$ in the region $0 < z < \De^{-\; 1/6}$ instead of $\pa^2 u/\pa x^2(x,t), 0 < x < t^{1/3}$.  We  have that
\be \label{R4}
g'(z) = \frac{-2\beta(z+\mu \sqrt{\De})}{\sig^2} \exp \left[ - \frac{(z + \mu\sqrt{\De})^2}{2\sig^2} \right]
\ee
\[    +\  G(z,\De) - \sqrt{2\pi \sig^2} \ e^{\rho \De} \ \De \; V''_0(z\sqrt{\De} - u_2 + u_0).	\]

We can estimate the function $G(z,\De)$ in the same way as we obtained the estimate (\ref{G4}).  In fact differentiating (\ref{F4}) and estimating as before we have the inequality,
\be \label{O4}
2 \left[ \frac x{\sig^2 t} + \frac{2\mu}{\sig^2}\right] \exp \left[ - \; \frac {(x + \mu t)^2}{2\sig^2 t} \right] \le \sqrt{2\pi\sig^2 t}  \ \frac{\pa^2 p}{\pa x^2} \; (x,t)
\ee
\[	\le 2 \left[ \frac x{\sig^2 t} + \frac{4\mu}{\sig^2}\right] \exp \left[ - \ \frac {(x + \mu t)^2}{2\sig^2 t}\right], \  \  x \ge 2\mu t.  \]
From (\ref{O4}) we can also estimate $\pa^2 h/\pa x^2$ using the identity,
\be \label{P4}
- \; \frac{\pa^2 h}{\pa x^2} \; (x,t) =\mu \int^t_0 \; \frac{\pa^2 p}{\pa x^2} \; (x,s)ds.
\ee
We can estimate the RHS of (\ref{P4}) by substituting the RHS of (\ref{O4}) and making the change of variable $w=(x + \mu s)/\sqrt{s}$.  Thus we have
\begin{eqnarray} \label{Q4}
- \; \frac{\pa^2 h}{\pa x^2} \; (x,t) &\le& \frac{32\mu}{\sig^2\sqrt{2\pi\sig^2}} \int^\infty_{(x+\mu\sqrt{t})/t} \exp \left[ - \; \frac{w^2}{2\sig^2} \right] dw \\
&\le& \frac{11}{\sqrt{2\pi\sig^2}} \exp \left[ - \; \frac{(x+\mu t)^2}{2\sig^2 t}\right], \ x \ge 2\mu t. \nn
\end{eqnarray}
The function $G(z,\De)$ is therefore bounded from the estimates (\ref{O4}), (\ref{Q4}) as
\be \label{S4}
0 \le -G(z,\De) \le \left( \frac{6\beta \mu}{\sig^2} + 11 \right)\sqrt{\De} \ \exp \left[ - \; \frac{(z+\mu \sqrt{\De})^2}{2\sig^2 }\right], \ z \ge 2\mu \sqrt{\De}.
\ee
From (\ref{R4}), (\ref{S4}) we see that there exists $\De_0 > 0$ such that for $0 < \De < \De_0$ one has $g'(z) < 0$ for  $0 < z < \sig\sqrt{2|\ln \; \De|}$ and $g'(z) > 0$ for $\sig\sqrt{3|\ln \; \De|} < z < \De^{-1/6}$.  Here there is a solution $z = a(\De)/\sqrt{\De}$ of the equation $g'(z) = 0$ in the interval $\sig\sqrt{2\ln \; \De} < z < \sig\sqrt{3\ln \; \De}$.  Evidently then  $\displaystyle{\lim_{\De \ra 0}} \; a(\De) = 0$.

We complete the proof by showing that  $g''(z) > 0$ for $\sig\sqrt{2|\ln \; \De|} < z < \sig \sqrt{3|\ln \; \De|}$, provided  $0 < \De < \De_0$.  This is accomplished by estimating $\pa^3p/\pa x^3$ in a similar way to how we estimated $\pa^2 p/\pa x^2$. 
\end{proof}
\begin{lem}  Let $u(x,t)$ be the function (\ref{K4}), and define $T > 0$ as $T = \sup\{ t>0 : \displaystyle{\sup_{x>0}} \;  u(x,t) > u(0,t)\}$.  Then for $0 < t < T$ there exists unique $x(t) > 0$ satisfying $\displaystyle {\lim_{t\ra 0}}\; x(t) = 0$, with the property that $\pa u/\pa x(x,t) > 0$ for $0 \le x < x(t)$ and $\pa u/\pa x (x,t) < 0$ for $x > x(t)$.  Furthermore $\pa^2 u/\pa x^2(x(t), t) < 0$.
\end{lem}
\begin{proof}  By the definition of $T$ one must have $\pa u/\pa x(x,t) > 0$ for some $x$ when $0 < t < T$.  Since $\pa u/\pa x(x,t) <0$ for $x$ large it follows from Theorem A1 that a unique $x(t)$ exists for $0 < t < T$.  Now we apply the argument in Theorem A3, using Lemma 4.4 to conclude that $\pa^2 u/\pa x^2(x(t), t) < 0$.  
\end{proof}
\begin{proposition} The function $V(x)$ defined by (\ref{V2}) is $C^1$ but not $C^2$ at $x = u_1$.
\end{proposition}
\begin{proof}  Let $\beta$ be as in (\ref{V2}) and $u(x,t)$ be the function (\ref{K4}). Then $u(u_1, \De) = \pa u/ \pa x(u_1, \De) = 0$, $u(0,\De) < 0, u(x,\De) \ra -\infty$ as $x \ra \infty$.  Then by Lemma 4.5 we have $\pa^2 u/\pa x^2 (u_1,\De) < 0$, whence $V$ is not $C^2$. 
\end{proof}

\appendix
\section{Some Consequences of the Maximum principle}

Here we prove some general results for the heat equation which are used in earlier sections.  Let $a(x,t), b(x,t)$ be uniformly bounded smooth functions in $(x,t), \; x \in \R, \; t \ge 0$, with the property that $a$ is also uniformly bounded from below by a positive constant.  We define the operator $L$ on $C^2$ functions $u(x,t)$ by
\[	Lu(x,t) = a(x,t) \ \frac{\pa^2u}{\pa x^2} + b(x,t) \frac{\pa u}{\pa x} - \frac{\pa u}{\pa t} .	\]
\begin{theorem} Suppose $u(x,t)$ is a $C^2$ function in $\{ (x,t) \in \R^2 : x \ge 0, \; t\ge 0\}$ satisfying $Lu \equiv 0$.  Suppose further that $u(0,t) = 0, \ t \ge 0$, and the set $\{ x > 0 : u(x,0) >0 \}$ is a semi-infinite interval.  Then for any $t > 0$ there is at most one point $x(t) > 0$ satisfying $u(x(t), \; t) = 0$.
\end{theorem}
\begin{proof}  Observe that since the function $u(x,0), \; x \ge 0$, can be negative only on a bounded set it follows that 
$u(x,t), \; x \ge 0$, is uniformly bounded below for all $t \ge 0$.  Further, one has $\displaystyle{\liminf_{x\ra \infty}} \; 
u(x,t) \ge 0$, for all $t \ge 0$. For some $T>0$ suppose there is an $x_0>0$ with $u(x_0, T) < 0$.  Let $D$ be the maximal domain containing $(x_0, T)$ such that 
$u(x,t) < 0$ for $(x,t) \in D$.  Evidently $u(x,t) = 0 $ for $(x,t) \in \pa D \cap \R\times(0,\infty)$.  Hence by the maximum principle the minimum of $u$ occurs at an interior point unless there exists $(x_1,0) \in \pa D$ with $u(x_1, 0) < 0$.  We conclude that there is a path $\ga(s), \ 0 \le s \le T$, with $\ga(T) = x_0$ and $u(\ga(s),s) < 0, \ 0 \le s \le T$.

Now we define $\Om$ to be the domain $\Om = \{(x,s) : 0 < x < \ga(s), \; 0 < s < T\}$.  Then $u \le 0$ on $\pa\Om \cap \R \times [0, T)$ and strictly negative on part of this boundary.  The maximum principle therefore implies that $u$ is strictly negative on $\pa\Om \cap \R^+ \times \{T\}$.  Thus $u$ is strictly negative on the interval $(0, x_0]$.  We define now $x(t)$ as $x(t)=0$ if $u(x,t) \ge 0, \ x > 0; \ x(t) = \infty$ if $u(x,t) < 0, x > 0; \ x(t)$ $= \lim\sup \{x>0 : u(x,t) <0\}$, otherwise.  If $0 < x(t) < \infty$ then the maximum principle implies that $u(x,t) > 0$ for $x > x(t)$. 
\end{proof}
\begin{theorem}  Suppose $u(x,t)$ is a $C^2$ function in $\{ (x,t) \in \R^2 : x \ge 0,  \ t \ge 0\}$ satisfying $Lu\equiv 0$.  Suppose further that $u(0,t) < 0, \; t \ge 0$, and that the set $\{ x > 0 : u(x,0) > 0\}$ is an open interval.  Then for any $t > 0$ the set $\{ x > 0 : u(x,t) \ge 0\}$ is either empty, a single point or a closed interval.  If the set is a closed interval then $u$ is strictly positive on its interior.
\end{theorem}
\begin{proof} By Theorem A1 we may assume that the set $\{ x > 0 : u(x,0) > 0\}$ is a finite interval.  Hence for all $t \ge 0$, $\displaystyle{\limsup_{x\ra \infty}}\ u(x,t) \le 0$.  Now let $D = \{ (x,t) :x > 0 ,  \ t >0,  \ u(x,t) > 0\}$.  Since $u$ is uniformly bounded above it follows by the maximum principle that $D$ is connected and that $\{(x,t) : x >0,  \ u(x,0) > 0\} \subset \pa D$.  We show that the intersection of $D$ with any line $t =$ constant is either empty or an open interval.  To do this let $a(t), b(t)$ be defined by
\[	a(t) = \inf \{ x > 0 : (x,t) \in D \},   \]
\[      b(t) = \sup \{x > 0: (x,t) \in D \}.	\]
We define a domain $\Om = \{(x,t) : t>0, \  a(t) < x < b(t) \}$, whence $D \subset \Om$.  For $T > 0$ let $D_T$ and $\Om_T$ be defined by $D_T = D \cap [\R \times (0, T)]$, $\Om_T = \Om \cap [\R \times (0, T)]$.  Suppose for some $T$ one has $D_T \not= \Om_T$.  By the maximum principle if the minimum of $u$ on $\Om_T$ is negative then it must be at a point $(x_0, T)$ with $a(T) < x_0 < b(T)$.  This contradicts Lemma 3 of Chapter 3, $\S2$ of \cite{pw}.  Hence $u \ge 0$ on $\Om_T$ and consequently by the maximum principle again $u > 0$ on $\Om_T$.  Applying the maximum principle to the complement of $\Om_T$ in $\R \times (0, T)$ we conclude that the set $\{x \in \R : u(x,t) \ge 0\} = [a(t), b(t)]$.

Finally suppose there is a minimum $T$ such that $\Om = \Om_T$.  As before the maximum principle implies that $u$ is strictly negative on the complement of $\bar \Om_T$ in $\R \times (0, \infty)$.  This completes the proof. 
\end{proof}
Next we wish to show that the solution $x(t)$ in Theorem A1 of the equation $u(x,t) = 0$ is nondegenerate.  As a consequence it follows that $x(t)$ is a smooth function of $t$ for $t > 0$.
\begin{theorem}  Suppose $u(x,t), \ x \ge 0, \ t \ge 0$, satisfies the assumptions of Theorem A1.  Let $W_0 = \{ x>0 : \pa u / \pa x (x,0) > 0 \}$ and assume that the boundary $\pa W_0$ of $W_0$ has no finite limit points.  Then there is the inequality $\pa u / \pa x (x(t), t) > 0, \; t > 0$.
\end{theorem}

\begin{proof}  For $t > 0$ let $W_t = \{ x>0 : \pa u / \pa x (x,t) > 0 \}$.  Since $\pa u / \pa x(x,t)$ also satisfies a diffusion equation it follows that $\pa W_t$ has no finite limit points.  Furthermore, for $t$ small and $x \not\in \pa W_t$ the derivative $\pa u / \pa x (x,t)$ is nonzero.

We first show that $\pa u / \pa x (x(t), t) > 0$ for $t >0$ sufficiently small.  We argue by contradiction.  Suppose $t_n, \; n=1,2,... $ is a positive sequence with $\displaystyle{\lim_{n\ra \infty}}\; t_n = 0$ such that $x(t_n) \in \pa W_{t_n}$, $n=1,2,\dots$.  We may assume wlog that $\displaystyle{\lim_{n\ra \infty}}\; x(t_n) = x_\infty$.  There is then a possibly nontrivial closed interval $[a, b]$ such that $x_\infty \in [a,b], \ \pa u/\pa x (x,0) = 0, \ x \in [a,b]$ and $\pa u/\pa x (x,0) \not= 0$, $x \not\in [a,b]$ but sufficiently close to $a$ or $b$.  Suppose now that $\pa u/\pa x(x,0) > 0$ for $x$ close to $a$ or close to $b$.  Then by the maximum principle it follows that $\pa u / \pa x (x,t) > 0$ for $x$ in a neighborhood of $[a,b]$ when $t$ is sufficiently small.  In particular $\pa u / \pa x (x(t_n), t_n) > 0$ for $n$ large, which is a contradiction.  Alternatively we can have $\pa u / \pa x (x,0) > 0$ for $x$ close to $a$ and $\pa u / \pa x (x,0) < 0$ for $x$ close to $b$.  In that case the maximum principle implies $u(x,t) < 0$ for $x$ in a neighborhood of $[a,b]$ when $t$ is small.  This again contradicts the fact that $u(x(t_n), t_n) = 0$ as $n \ra \infty$.  We conclude therefore that $\pa u / \pa x (x(t), t) > 0$ for $t > 0$ sufficiently small.

Finally we show that $\pa u / \pa x (x(t), t) > 0$ for all $t$.  To see this consider $T > 0$ and suppose that $\pa u / \pa x (x(t), t) > 0$ for all $0 < t < T$, $\displaystyle{\lim_{t \ra T}}\; x(t) = x_T$.  If $x_T = \infty$ then $u(x,T) \le 0, x \ge 0$.  If $x_T = 0$ then $u(x,T) \ge 0, x \ge 0$.  Suppose now $0 < x_T < \infty$, whence $u(x,T), \; x > 0$, takes on both positive and negative values.  For $0 < t < T$ there exists $a(t), b(t)$ with $0 \le a(t) < x(t) < b(t) \le \infty$ such that $\pa u / \pa x (x, t) > 0$ for $a(t) < x < b(t)$ and $a(t), b(t) \in \pa W_t$ if $a(t) > 0$ and $b(t) < \infty$.  Since $u(0,t) = 0$ it follows that $\pa u / \pa x (x, t) < 0$ for $x > 0$ small.  Hence $a(t) > 0$.

We show that $x_T$ cannot be the limit as $t \ra T$ of any points in $\pa W_t \cap [0,x(t)]$. Suppose first that $\pa W_t \cap [0,x(t)] = \{a(t)\}$.  Then $\pa u / \pa x(x,t), \; 0 < x < x(t)$, has just one sign change, whence $\pa u / \pa x(x,T), \; 0 < x < x_T$,  has also at most one sign change.  Since $u(x_T, T) = u(0,T) = 0$ it follows that $\displaystyle{\lim_{t \ra T}} \ a(t) < x_T$.  Alternatively let $a^*(t) \in \pa W_t \cap [0,x(t)]$ have the property that
\[	u(a^*(t), t) = \sup \{ u(x,t) : x \in \pa W_t \cap [0, x(t)] \} .  \]
Evidently $u(a^* (t),t) < 0$.  By the maximum principle it follows that $u(a^* (t),t)$ is a decreasing function of $t$.  Hence $x_T$ cannot be the limit as $t \ra T$ of any point in $\pa W_t \cap [0,x(t)]$.

Since we can argue similarly that $x_T$ is not the limit as $t \ra T$ of any point in $\pa W_t \cap [x(t), \infty)$ we conclude that $\pa u / \pa x (x(T), T) > 0$.  We have shown that the set $\{ t > 0 : \pa u / \pa x(x(t), t) > 0 \}$ is both open and closed whence it must be the interval $(0, \infty)$. 
\end{proof}

\thanks{ {\bf Acknowledgement:} The authors would like to thank Jeffrey Rauch and Peter Bates for helpful conversations. This research was 
partially supported by NSF under grants DMS-0138519 and DMS-0500608.


\begin{thebibliography}{9}
\bibitem{b}
M. Bramson,
\newblock \textit{Convergence of solutions of the Kolmogorov equation to travelling waves},
\newblock Mem. Amer. Math. Soc. \textbf{44}:285 (1983).
\bibitem{fr} W. Fleming and R. Rishel,
\newblock Deterministic and stochastic optimal control.
\newblock Springer Verlag, New York 1975, 222 pp.
\bibitem{fls} W. Fleming and H. Soner,
\newblock  Controlled Markov Processes and Viscosity Solutions.
\newblock  Springer Verlag, New York 1993, 428 pp.
\bibitem{fs} A. Friedman and W. Shen,
\newblock \textit{A variational inequality approach to financial valuation of retirement benefits based on salary}
\newblock Finance and Stochastics \textbf{6} (2002), 273-302.
\bibitem{ht} B. Hojgaard and M. Taksar,
\newblock \textit{Controlling risk exposure and dividends payout schemes: insurance company example},
\newblock Math. Finance \textbf{9} (1999), 153-182.
\bibitem{ht1} B. Hojgaard and M. Taksar,
\newblock \textit{Optimal risk control for a large corporation in the presence of returns on investments},
\newblock Finance and Stochastics \textbf{5} (2001), 527-547.
\bibitem{mck}
H. McKean, 
\newblock \textit{Application of Brownian Motion to the equation of Kolmogorov-Petrowskii-Piskounov},
\newblock Comm. Pure Appl. Math. \textbf{28} (1975), 323-331.
\bibitem{mr}
A. Milne and D. Robertson,
\newblock \textit{Firm behavior under the threat of liquidation},
\newblock J. Economic Dynamics and Control \textbf{20} (1996), 1427-1449.
\bibitem{pw} M. Protter and H. Weinberger,
\newblock   Maximum principle in differential equations. Corrected reprint of the 1967 original.
\newblock Springer Verlag, New York 1984, 261 pp.
\bibitem{pk} S. Peura and J. Keppo,
\newblock  \textit{Optimal bank capital with costly recapitalization},
\newblock 2005 preprint.
\end{thebibliography}
\end{document}